\documentclass[11pt,leqno]{amsart}

\newtheorem{lemma}{Lemma}
\newtheorem{theorem}{Theorem}
\newtheorem{corollary}{Corollary}
\newtheorem{proposition}{Proposition}
\newtheorem{definition}{Definition}

\newtheorem{remark}{Remark}
\usepackage{graphics}
\usepackage{epsfig}
\usepackage{amssymb,amsfonts,amscd}
\usepackage{enumerate}
\def\B{{\mathbb B}}
\def\C{{\mathbb C}}

\def\R{{\mathbb R}}
\def\CC{{\mathcal C}}
\begin{document}
                      
\title[Riemann maps in almost complex manifolds]
{Riemann maps in almost complex manifolds}
\author{B.Coupet, H.Gaussier and A.Sukhov}

\address{\begin{tabular}{lll}
Bernard Coupet & & Herv\'e Gaussier\\
LATP UMR 6632 & & LATP UMR 6632\\
C.M.I. & & C.M.I.\\
39, rue Joliot-Curie & &39, rue Joliot-Curie\\
13453 Marseille Cedex 13 & &13453 Marseille Cedex 13\\
FRANCE & & FRANCE\\
{\rm coupet@cmi.univ-mrs.fr} & & {\rm gaussier@cmi.univ-mrs.fr}\\
 & & \\
Alexandre Sukhov & & \\
AGAT UMR 8524 & & \\
U.S.T.L.& & \\
Cit\'e Scientifique & & \\
59655 Villeneuve d'Ascq Cedex & & \\
FRANCE & & \\
{\rm sukhov@agat.univ-lille1.fr} & &
\end{tabular}
}

\email{} 


\date{\number\year-\number\month-\number\day}

\begin{abstract}
We prove the existence of stationary discs in the ball for small 
almost complex deformations of the standard structure. We define a local
analogue of the 
Riemann map and establish its main properties. 
These constructions are applied to study the local geometry of
almost complex manifolds and their morphisms. 
\end{abstract}

\maketitle

\section*{Introduction}
The notion of stationary disc was introduced by L.Lempert
in~\cite{le81}. A holomorphic disc in a domain $\Omega$ in
$\C^n$ is a holomorphic map from the unit disc $\Delta$ to $\Omega$.  A
proper holomorphic disc $f: \Delta \rightarrow \Omega$, continuous up
to $\partial \Delta$, is called {\sl stationary} if there is a
meromorphic lift $\hat{f}$ of $f$ to the cotangent bundle $T^*(\C^n)$
with an only possible pole of order at most one at origin, such that
the image $\hat{f}(\partial \Delta)$ is contained in the conormal
bundle $\Sigma(\partial \Omega)$ of $\partial \Omega$. 
L.Lempert proved in \cite{le81} that for a
strictly convex domain stationary discs coincide with extremal discs
for the Kobayashi metric and studied their basic properties. 
Using these discs he introduced a multi dimensional analogue of the
Riemann map (see also \cite{se92} by S.Semmes and \cite{ba-le98} by
Z. Balogh and Ch. Leuenberger). The importance of this object comes
from its links with the
complex potential theory \cite{le81,le83}, moduli spaces of Cauchy-Riemann
structures and contact geometry \cite{bl94,bl-du91,bl-du-ka87,le88,se92},
and the mapping problem \cite{le86,le91,tu01}.

Stationary discs are natural global biholomorphic invariants of complex
manifolds with boundary. In view of the recent progress in symplectic geometry
due to the application of almost complex structures and pseudo holomorphic
curves, it seems relevant to find an analogue of 
Lempert's theory in the almost complex
case. In this paper we restrict ourselves to small almost complex deformations
of the standard integrable structure in $\C^n$. 
Since any almost complex structure
may be represented locally in such a form, our results can be viewed as
a local analogue of Lempert's theory in almost complex manifolds.
We prove the existence of stationary discs in the unit ball for small
deformations of the structure and show that they form a foliation of the 
unit ball (singular at the origin) in Section~4. 
Then we define an analogue of the Riemann map and establish
its main properties: regularity, holomorphy along the leaves 
and commutativity with (pseudo) biholomorphic maps (Theorem~\ref{TH1}).
As an application we prove the following criterion for the boundary regularity 
of diffeomorphisms in almost complex manifolds, which is a 
partial generalization of Fefferman's
theorem~\cite{fe74}: let $M$ and $M'$ be two $\mathcal C^{\infty}$ smooth real
$2n$-dimensional manifolds, $D \subset M$ and $D' \subset M'$ be
relatively compact domains. Suppose that there exists an almost
complex structure $J$ of class $\mathcal C^{\infty}$ on $\bar D$ such
that $(D,J)$ is strictly pseudoconvex. Then a $\mathcal C^{1}$
diffeomorphism $\phi$ between $\bar{D}$ and $\bar{D}'$ is of class 
$\mathcal C^{\infty}(\bar D)$ if and only if the
direct image $J':= \phi^*(J)$ is of class $\mathcal
C^{\infty}(\bar{D'})$ and $(D',J')$ is strictly pseudoconvex. (See Theorem~4).
This result is partially motivated by Eliashberg's question on the existence
of a symplectic analogue of Fefferman's theorem see~\cite{be90}. 
Finally we point out that
the Riemann map contains an essential information on the local geometry of an
almost complex manifold. It can be used to solve the local equivalence problem
for almost complex structures (Theorem~5). 
\vskip 0,1cm
\noindent{\it Acknowledgment.} The authors are grateful to Evgueni Chirka for
several helpful conversations.

\section{Preliminaries}

Everywhere in this paper $\Delta$ denotes the unit disc in $\C$ and $\B_n$
the unit ball in $\C^n$. Let
$(M,J)$ be an almost complex manifold ($J$ is a smooth $\mathcal
C^\infty$-field on the tangent bundle $TM$ of $M$, satisfying
$J^2=-I$).  By an abuse of notation $J_0$ is the standard structure on
$\C^k$ for every integer $k$.  A $J$-holomorphic disc in $M$ is a
smooth map from $\Delta$ to $M$, satisfying the quasilinear elliptic
equation $df \circ J_0 = J \circ df$ on $\Delta$.  

 An important special case of an almost complex manifold is a bounded
 domain $D$ in $\C^n$ equipped with an almost complex structure $J$,
 defined in a neighborhood of $\bar{D}$, and sufficiently close to the
 standard structure $J_0$ in the $\mathcal C^2$ norm on
 $\bar{D}$. Every almost complex manifold may be represented locally
 in such a form. More precisely, we have the following Lemma. 

\begin{lemma}
\label{suplem1}
Let $(M,J)$ be an almost complex manifold. Then for every point $p \in
M$ and every $\lambda_0 > 0$ there exist a neighborhood $U$ of $p$ and a
coordinate diffeomorphism $z: U \rightarrow \mathbb B$ such that
$z(p) = 0$, $dz(p) \circ J(p) \circ dz^{-1}(0) = J_0$  and the
direct image $J' = z^*(J)$ satisfies $\vert\vert J' - J_0
\vert\vert_{\CC^2(\overline {\mathbb B})} \leq \lambda_0$.
\end{lemma}
\proof There exists a diffeomorphism $z$ from a neighborhood $U'$ of
$p \in M$ onto $\mathbb B$ satisfying $z(p) = 0$ and $dz(p) \circ J(p)
\circ dz^{-1}(0) = J_0$. For $\lambda > 0$ consider the dilation
$d_{\lambda}: t \mapsto \lambda^{-1}t$ in $\C^n$ and the composition
$z_{\lambda} = d_{\lambda} \circ z$. Then $\lim_{\lambda \rightarrow
0} \vert\vert z_{\lambda}^{*}(J) - J_0 \vert\vert_{\CC^2(\overline
{\mathbb B})} = 0$. Setting $U = z^{-1}_{\lambda}(\mathbb B)$ for
$\lambda > 0$ small enough, we obtain the desired statement. \qed

In particular, every almost complex structure $J$ sufficiently close 
to the standard
structure $J_0$ will be written locally 
$J=J_0 + \mathcal O(\|z\|)$. 
Finally by a small perturbation 
(or deformation) of the
standard structure $J_0$ defined in a neighborhood of $\bar{D}$, where
$D$ is a domain in $\C^n$, we will mean a smooth one parameter family
$(J_\lambda)_\lambda$ of almost complex structures defined in a
neighborhood of $\bar{D}$, the real parameter $\lambda$ belonging to a
neighborhood of the origin, and satisfying~: $\lim_{\lambda
\rightarrow 0}\|J_\lambda - J_0\|_{\CC^2(\bar{D})} = 0$.

\vskip 0,1cm
Let $(M,J)$ be an almost complex manifold. We denote by $T(M)$ the real 
tangent bundle of $M$ and by $T_\C(M)$ its complexification. 
Then~:
$$
T_{(1,0)}(M):=\{ X \in T_\C (M) : JX=iX\} = 
\{\zeta -iJ \zeta, \zeta \in T(M)\},
$$
$$
T_{(0,1)}(M):=\{ X \in T_\C (M) : JX=-iX\} = 
\{\zeta +iJ \zeta, \zeta \in T(M)\}.
$$
Let $T^*(M)$ denote the cotangent bundle of $M$. Identifying $\C
\otimes T^*M$ with $T_\C^*(M):=Hom(T_\C M,\C)$ we may consider a complex
one form on $M$ as an element in $T_\C^* (M)$. 
We define the set of complex forms of type $(1,0)$ on $M$ by~:
$$
T^*_{(1,0)}(M)=\{w \in T_\C^* (M) : w(X) = 0, \forall X \in T_{(0,1)}(M)\}
$$
and the set of complex forms of type $(0,1)$ on $M$ by~:
$$
T^*_{(0,1)}(M)=\{w \in T_\C^* (M) : w(X) = 0, \forall X \in T_{(1,0)}(M)\}.
$$
Let $\Gamma$ be a real smooth submanifold in $M$ and let $p \in
\Gamma$. We denote by $H^J(\Gamma)$ the $J$-holomorphic tangent bundle 
$T\Gamma \cap JT\Gamma$.

\begin{definition}
A real submanifold $\Gamma$ in $M$ is {\sl totally real} if 
$H_p^J(\Gamma)=\{0\}$ for every $p \in \Gamma$.
\end{definition}

In what follows we will need the notion of the Levi form of a hypersurface.

\begin{definition}\label{DEF}
Let $\Gamma=\{r=0\}$ be a smooth real hypersurface in $M$ 
($r$ is any smooth defining function of $\Gamma$) and let $p \in \Gamma$. 

$(i)$ The {\sl Levi form} of $\Gamma$ at $p$ is the map defined on
$H^J_p(\Gamma)$  by $\mathcal L_\Gamma^J(X_p) = (J^*dr)[\bar{X},X]_p$,
where  a vector field $X$ is a section of the $J$-holomorphic tangent
bundle  $H^J \Gamma$ such that $X(p) = X_p$.

$(ii)$ A real smooth hypersurface $\Gamma=\{r=0\}$ in $M$ is 
{\sl strictly $J$-pseudoconvex} if its Levi form $\mathcal L_\Gamma^J$ 
is positive definite on $H^J(\Gamma)$.
\end{definition}
 
We recall the notion of conormal bundle of a real submanifold in $\C^n$ 
(\cite{tu01}).
Let $\pi: T^*(\C^n) \rightarrow
\C^n$ be the natural projection. Then $T^*_{(1,0)}(\C^n)$ can be 
canonically identified with 
the cotangent bundle of $\C^n$.  
In the canonical complex coordinates $(z,t)$ on $T^*_{(1,0)}(\C^n)$
an element of the fiber at $z$ is a holomorphic form $w=\sum_j t_j dz^j$. 

Let $N$ be a real smooth generic submanifold in $\C^n$. 
The conormal bundle $\Sigma(N)$ of $N$ is
a real subbundle of $T^*_{(1,0)}(\C^n)|_N$ whose fiber at $z\in N$
is defined by $\Sigma_{z}(N) = \{ \phi \in T^*_{(1,0)}(\C^n): Re
\,\phi \vert T_{(1,0)}(N) = 0 \}$.

Let $\rho_1,\dots,\rho_d$ be local defining functions of $N$.
Then the forms $\partial \rho_1,\dots, \partial \rho_d$ form
a basis in $\Sigma_{z}(N)$ and every section $\phi$ of the bundle
$\Sigma(N)$ has the form $\phi = \sum_{j=1}^d c_j \partial \rho_j$, 
$c_1,\dots,c_d \in
\R$. We will use the following
(see\cite{tu01}):
\begin{lemma}\label{lemma}
Let $\Gamma$ be a real $\CC^2$ hypersurface in $\C^n$. 
The conormal bundle $\Sigma(\Gamma)$ is 
a totally real
submanifold of dimension $2n$ in $T^*_{(1,0)}(\C^n)$ if and 
only if the Levi form of $\Gamma$ is nondegenerate.
\end{lemma} 
In Section 5 we will introduce an analogue of this notion in the
almost complex case.

\section{Existence of discs attached to a  real submanifold 
of an almost complex manifold}

\subsection{ Partial indices and the Riemann-Hilbert problem}

In this section we introduce basic tools of the linear Riemann-Hilbert 
problem.

Let $V \subset \C^N$ be an open set. We denote
by $\CC^k(V)$ the Banach space of (real or complex valued) functions
of class $\CC^k$ on $V$ with the standard norm
$$
\parallel r \parallel_k =
\sum_{\vert \nu \vert \leq k}
\sup \{ \vert D^{\nu} r(w) \vert : w \in V \}.
$$
For a positive real 
number $\alpha <1$ and a Banach space $X$, we denote by
$\CC^{\alpha}(\partial \Delta,X)$ the Banach space of all functions
$f: \partial \Delta \rightarrow X$ such that

$$
\parallel f \parallel_{\alpha} :=
\sup_{\zeta \in \partial \Delta} \parallel f(\zeta) \parallel +
\sup_{\theta,\eta \in \partial \Delta, \theta \neq \eta}
\frac{\parallel f(\theta) - f(\eta)\parallel}{\vert \theta - \eta
\vert^{\alpha}} < \infty.
$$

\vskip 0,1cm \noindent If $\alpha = m +
\beta$ with an integer $m \geq 0$ and $\beta \in ]0,1[$, then we consider
the Banach space
$$
\CC^{\alpha}(V):=\{r \in \CC^m(V,\R): D^{\nu} r \in C^{\beta}(V), 
\nu: \vert \nu \vert \leq m\}
$$
and we set $\parallel r
\parallel_{\alpha} = 
\sum_{\vert \nu \vert} \parallel D^{\nu}r \parallel_{\beta}$.

Then a map $f$ is in $\CC^{\alpha}(V,\C^k)$ if and only if its components 
belong to $\CC^{\alpha}(V)$ and we say that $f$ is of class $\CC^{\alpha}$. 

\vskip 0,1cm
Consider the following situation:

\vskip 0,1cm
\noindent $\bullet$ $B$ is an open ball centered at the origin in $\C^N$ and 
$r^1,\dots,r^N$ are smooth $\CC^\infty$ functions defined in a neighborhood 
of  $\partial \Delta \times B$ in $\C^N \times \C$

\noindent $\bullet$ $f$ is a map of class $\CC^{\alpha}$ 
from $\partial \Delta$ to $B$, where $\alpha>1$ is a noninteger real
number

\noindent $\bullet$ for every $\zeta \in \partial \Delta$

\begin{itemize}
\item[(i)] $E(\zeta) = \{ z \in B: r^j(z,\zeta) = 0, 1 \leq j \leq N \}$ 
is a maximal totally real submanifold in $\C^N$,
\item[(ii)] $f(\zeta) \in E(\zeta)$,
\item[(iii)] $\partial_z r^1(z,\zeta) \wedge \cdots \wedge \partial_z
r^N(z,\zeta) \neq 0$ on $B \times \partial \Delta$.
\end{itemize}

Such a family ${E} = \{ E(\zeta) \}$
of manifolds with a fixed disc $f$ is called a
{\it totally real fibration} over the unit circle.
A disc attached to a fixed totally real manifold ($E$ is independent of 
$\zeta$) is a special case of a totally real fibration. 

\vskip 0,1cm
Assume that the defining function $r:=(r^1,\dots,r^N)$ of $E$ 
depends smoothly on a small real parameter $\varepsilon$, namely $r =
r(z,\zeta,\varepsilon)$, 
and that the fibration $E_0:=E(\zeta,0)$ corresponding
to $\varepsilon = 0$ coincides with the above fibration
$E$. Then for every sufficiently small $\varepsilon$ and for every
$\zeta \in \partial \Delta$ the manifold $E_\varepsilon:=E(\zeta,\varepsilon) 
:= \{ z \in B:
r(z,\zeta,\varepsilon) = 0\}$ is totally real. 
We call $E_{\varepsilon}$ a {\it smooth totally real deformation }  
of the totally
real fibration $E$. By a holomorphic disc $\tilde f$ attached to
$E_{\varepsilon}$  we mean a holomorphic map $\tilde f:\Delta
\rightarrow B$, continuous on $\bar\Delta$, satisfying
$r(f(\zeta),\zeta,\varepsilon) = 0$ on $\partial \Delta$. 

For every positive real noninteger $\alpha$ we denote by
$(\mathcal A^\alpha)^N$ the space of maps defined on
$\bar{\Delta}$, $J_0$-holomorphic on $\Delta$, and belonging to
$(\mathcal C^\alpha(\bar{\Delta}))^N$.

\subsection{Almost complex perturbation of discs}
We recall that for $\lambda$ small enough the
$(J_0,J_{\lambda} )$-holomorphy condition for a map $f:\Delta
\rightarrow \C^N$ may be written in the form

\begin{equation}\label{equa0}
\bar\partial_{J_{\lambda}} f = \bar\partial f +
q(\lambda,f)\partial f = 0
\end{equation}
where $q$ is a smooth function satisfying $q(0,\cdot) \equiv 0$, uniquely
determined by $J_\lambda$ (\cite{si94}).

Let $E_\varepsilon=\{r_j(z,\zeta,\varepsilon) = 0, 1 \leq j \leq
N\}$ be a smooth totally real deformation of a totally real fibration $E$. 
A disc $f \in \mathcal (\CC^\alpha(\bar{\Delta}))^N$ is attached
to $E_\varepsilon$ and is $J_\lambda$-holomorphic if and only if it
satisfies the following nonlinear boundary Riemann-Hilbert type  problem~:
$$
\left\{
\begin{array}{lll}
r(f(\zeta),\zeta,\varepsilon) &=& 0, \ \ \ \zeta \in \partial \Delta\\
 & & \\
\bar{\partial}_{J_\lambda}f(\zeta) &=& 0, \ \ \ \zeta \in \Delta.
\end{array}
\right.
$$
Let $f^0 \in \mathcal (\mathcal A^\alpha)^N$ be a disc attached to
$E$ and let $\mathcal U$ be a neighborhood of $(f^0,0,0)$ in the 
space $(\mathcal C^\alpha(\bar{\Delta}))^N \times \R \times \R$.
Given $(f,\varepsilon,\lambda)$ in $U$ define the maps 
$v_{f,\varepsilon,\lambda}: \zeta \in \partial 
\Delta \mapsto r(f(\zeta), \zeta, \varepsilon)$
and
$$
\begin{array}{llcll}
u &:& \mathcal U & \rightarrow & (\mathcal C^\alpha(\partial \Delta))^N \times
 \mathcal C^{\alpha-1}(\Delta)\\
  & & (f,\varepsilon,\lambda) & \mapsto & (v_{f,\varepsilon,\lambda}, 
\bar{\partial}_{J_\lambda}f).
\end{array}
$$

Denote by $X$ the Banach space $(\mathcal C^\alpha(\bar \Delta))^N$.
Since $r$ is of class $\CC^{\infty}$, 
the map
$u$ is smooth and the tangent map $D_Xu(f^0,0,0)$ (we consider 
the derivative
with respect to the space $X$) is a linear map from $X$ to 
$(\mathcal C^\alpha(\partial \Delta))^N \times \mathcal C^{\alpha-1}(\Delta)$,
defined for every $h \in X$ by 
$$\begin{array}{llllll}
D_Xu(f^0,0,0)(h) = \left(
\begin{matrix}
2 Re [G h] \\
\bar\partial_{J_0} h
\end{matrix}
\right),
\end{array}$$
where for $\zeta \in \partial \Delta$
$$\begin{array}{lllllll}
G(\zeta) = \left(
\begin{matrix}
\frac{\partial r_1}{\partial z^1}(f^0(\zeta),0)  &\cdots&\frac{\partial
r_1}{\partial z^N}(f^0(\zeta),0)\\
\cdots&\cdots&\cdots\\
\frac{\partial r_N}{\partial z^1}(f^0(\zeta),0)& \cdots&\frac{\partial r_N}
{\partial z^N}(f^0(\zeta),0)
\end{matrix}
\right)
\end{array}$$
(see \cite{gl94}).
\begin{proposition}\label{tthh}
Let $f^0:\bar \Delta \rightarrow \C^N$ be a $J_0$-holomorphic
disc attached to a totally real fibration $E$ in $\C^N$.  Let
$E_\varepsilon$ be a smooth totally real deformation 
of $E$ and
$J_\lambda$ be a smooth almost complex deformation of $J_0$ in a
neighborhood of $f(\bar{\Delta})$.
Assume that for some $\alpha > 1$
the linear map from $(\mathcal A^{\alpha})^N$ to $(\mathcal 
C^{\alpha-1}(\Delta))^N$ given by 
$h \mapsto 2 Re [G h]$
is surjective and has a $k$ dimensional kernel.
Then there exist $\delta_0,
\varepsilon_0, \lambda_0 >0$ such that for every $0 \leq \varepsilon
\leq \varepsilon_0$ and for every $0 \leq \lambda \leq \lambda_0$,
the set of $J_\lambda$-holomorphic discs $f$ attached to $E_\varepsilon$
and such that $\parallel f -f^0 \parallel_{\alpha} \leq \delta_0$ forms
a smooth $(N+k)$-dimensional 
submanifold 
$\mathcal A_{\varepsilon,\lambda}$ in the Banach space 
$(\mathcal C^\alpha(\bar{\Delta}))^N$.
\end{proposition} 

\noindent{\bf Proof.} According to the implicit function Theorem, 
the proof of Proposition~\ref{tthh} reduces to the proof of the 
surjectivity of $D_Xu$. 
It follows by classical
one-variable results on the resolution of the
$\bar\partial$-problem in the unit disc that the linear map from
$X$ to $\mathcal C^{\alpha-1}(\Delta)$ given by 
$h \mapsto \bar \partial h$
is surjective. More precisely, given $g \in \mathcal C^{\alpha-1}(\Delta)$ 
consider
the Cauchy transform 

$$T_{\Delta}(g) : \tau \in \partial \Delta \mapsto   
\int\int_{\Delta} \frac{g(\zeta)}{\zeta - \tau}d\zeta d\bar{\zeta}.$$

For every function $g \in \mathcal
C^{\alpha-1}(\Delta)$ the solutions $h \in X$ of the equation
$\bar\partial h = g$ have the form $h = h_0 + T_{\Delta}(g)$
where $h_0$  is an arbitrary function in $({\mathcal A}^{\alpha})^N$. 
Consider the equation

\begin{equation}\label{equa1}
D_Xu(f^0,0,0)(h) = \left(
\begin{matrix}
g_1 \\
g_2
\end{matrix}
\right),
\end{equation}
where $(g_1,g_2)$ is a vector-valued function with components 
$g_1 \in \mathcal C^{\alpha-1}(\partial \Delta)$ and $g_2 \in \mathcal 
C^{\alpha-1}(\Delta)$. Solving the
$\bar\partial$-equation for the second component, we reduce 
(\ref{equa1}) to 
$$
2 Re [G(\zeta) h_0(\zeta)] = g_1 - 2 Re [G(\zeta) T_{\Delta}(g_2)(\zeta)]
$$
with respect to $h_0 \in (\mathcal A^{\alpha})^N$. 
The surjectivity of the map $ h_0 \mapsto 2 Re [G h_0]$ gives the result. \qed

\subsection{Riemann-Hilbert problem on the (co)tangent bundle of an almost 
complex manifold} 
Let $(J_{\lambda})_\lambda$ be an almost complex deformation 
of the standard structure $J_0$ on $\B_n$, satisfying $J_\lambda(0)= J_0$ and 
consider a (1,1) tensor field defined  on the bundle $\mathbb B_n
\times \R^{2n}$, represented in the $(x,y)$ coordinates by a $(4n \times
4n)$-matrix

$$
T_{\lambda} = \left(
\begin{matrix}
J_{\lambda}(x) & 0\\
\sum y_k A^k_{\lambda}(x)& B_{\lambda}(x)
\end{matrix}
\right),
$$
where $A_{\lambda}^k(x) = A^k(\lambda,x)$, $B_{\lambda}(x) =
B^k(\lambda,x)$ are smooth $(2n \times 2n)$-matrix functions 
($B_{\lambda}$ will be either $J_{\lambda}$ or ${}^{t}J_{\lambda}$). We
stress that we do not assume $T_{\lambda}$ to be an almost complex structure,
namely we do not require the identity $T_{\lambda}^2 = -Id$. 
In what follows we always assume that for every $k$: 

\begin{equation}
\label{norm1}
A^k_0(x) \equiv 0, \ \ \ {\rm for \ every}\ k
\end{equation}
and that one of the following two conditions holds:

\begin{equation}\label{norm2}
B_0(x)  = B_{\lambda}(0)  =  J_0, \ {\rm for \ every}\ \lambda,\ x,
\end{equation}
\begin{equation}\label{norm2'}
B_0(x)  =  B_{\lambda}(0)  =  ^tJ_0, \ {\rm for \ every}\ \lambda,\ x. 
\end{equation}
In what follows the trivial bundle $\mathbb B \times \R^{2n}$
over the unit ball will be a local coordinate
representation of the tangent or cotangent bundle of an almost complex
manifold. We denote by $x = (x^1,\dots,x^{2n})\in \mathbb B_n$
and $y = (y_1,\dots,y_{2n}) \in \R^{2n}$ the coordinates on the base and
fibers respectively. We identify the base space $(\R^{2n},x)$
with $(\C^n,z)$. 
Since ${}^tJ_0$ is orthogonally equivalent to $J_0$ we may identify
$(\mathbb R^{2n},{}^tJ_0)$ with $(\mathbb C^n,J_0)$. After this identification
the ${{}^tJ_0}$-holomorphy is expressed by the
$\bar{\partial}$-equation in the usual $t$ coordinates in $\C^n$.

By analogy with the almost complex case we say that a smooth map
$\hat{f}=(f,g): \Delta \rightarrow \mathbb B \times \R^{2n}$ is 
{\it $(J_0,T_{\lambda})$-holomorphic} if it satisfies 

$$
T_{\lambda}(f,g) \circ  d\hat{f} = d\hat{f} \circ J_0
$$
on $\Delta$.

For $\lambda$ small enough this can be rewritten as the
following Beltrami type quasilinear elliptic equation:

$$
(\mathcal E)\left\{
\begin{array}{cll}
\bar \partial f + q_1(\lambda,f))\partial f & = & 0\\
 & & \\
\bar \partial g + q_2(\lambda,f))\partial g + q_3(\lambda,f) g
 & = & 0,
\end{array}
\right.
$$
where the first equation coincides with the $(J_0,J_{\lambda})$-holomorphy
condition for $f$ that is $\bar\partial_{J_\lambda} f 
= \bar \partial f + q_1(\lambda,f))\partial$.
The coefficient $q_1$ is uniquely determined by $J_{\lambda}$ and, in view 
of (\ref{norm1}), (\ref{norm2}), (\ref{norm2'}) 
the coefficient $q_k$ satisfies, for $k=2,3$:

\begin{equation}
\label{nprm3}
q_k(0,\cdot) \equiv 0, \ q_k(\cdot,0) \equiv 0.
\end{equation}

We point out that in $(\mathcal E)$ 
the equations for the fiber component $g$ are
obtained as a small perturbation of the $\bar\partial$-operator. 
An important feature of this system is that the second equation is
{\it linear} with respect to the fiber component $g$. 

\begin{definition}
We call the above tensor field $T_{\lambda}$ a prolongation of the
structure $J_{\lambda}$ to the bundle $\mathbb B_n \times \R^{2n}$. The
operator

$$\bar\partial_{T_{\lambda}} :     \left(
\begin{matrix}
f\\
g
\end{matrix}
\right) \mapsto  \left(
\begin{matrix}
\bar \partial f + q_1(\lambda,f))\partial f\\
\bar \partial g + q_2(\lambda,f))\partial g + q_3(\lambda,f) g
\end{matrix}
\right)$$
is called an elliptic prolongation of
the operator $\bar\partial_{J_\lambda}$ to the bundle $\mathbb B_n \times
\R^{2n}$ associated with $T_{\lambda}$.
\end{definition}
 
Let $r^j(z,t,\lambda)$ , $j = 1,\dots,4n$ be $\CC^{\infty}$-smooth real 
functions
on $\mathbb B \times \mathbb B \times [0,\lambda_0]$ and let 
$r:=(r^1,\dots,r^{4n})$. 
In what follows we consider the following nonlinear boundary
Riemann-Hilbert type problem for the operator
$\bar\partial_{T_{\lambda}}$:
$$
({\mathcal
BP}_{\lambda})
\left\{
\begin{array}{cll}
r(f(\zeta),\zeta^{-1}g(\zeta),\lambda) &=& 0 \ \ {\rm on} \ \partial \Delta,\\
 \\
\bar \partial_{T_\lambda} (f,g) &=& 0,
\end{array}
\right.
$$
on the space $\CC^{\alpha}(\bar\Delta,\B_n\times \B_n)$.

The boundary problem $({\mathcal
BP}_{\lambda})$ has the following geometric sense.
Consider the disc $(\hat f,\hat g) := (f, \zeta^{-1}g)$ on 
$\Delta\backslash \{0\}$ and the set
$E_{\lambda} := \{ (z,t): r(z,t,\lambda) = 0 \}$. The
boundary condition in $({\mathcal BP}_{\lambda})$ means that 
$$
(\hat f, \hat g)(\partial \Delta) \subset E_{\lambda}.
$$
This boundary problem has the following {\it invariance
property}. Let $(f,g)$ be a solution of
$({\mathcal BP}_{\lambda})$ and let $\phi$ be a automorphism of 
$\Delta$. 
Then $(f \circ \phi,c g \circ \phi)$ also satisfies the
$\bar\partial_{T_{\lambda}}$ equation for every complex constant $c$.
In particular, if $\theta \in[0,2\pi]$ is fixed, then the disc
$(f(e^{i\theta}\zeta),e^{-i\theta}g(e^{i\theta}\zeta))$ satisfies the
$\bar\partial_{T_{\lambda}}$-equation on $\Delta \backslash \{0\}$ and
the boundary of the disc
$(f(e^{i\theta}\zeta),e^{-i\theta}\zeta^{-1}g(e^{i\theta}\zeta))$ is
attached to $E_{\lambda}$.  This implies the following 

\begin{lemma}
\label{circle}
If $(f,g)$ is a solution of $({\mathcal BP}_{\lambda})$,
then 
$\zeta \mapsto (f(e^{i\theta}\zeta),e^{-i\theta}g(e^{i\theta}\zeta))$ is also a
solution of $({\mathcal BP}_{\lambda})$.
\end{lemma}

Suppose that this problem has a solution $(f^0,g^0)$ for $\lambda = 0$
(in view of the above assumptions this solution is holomorphic on
$\Delta$ with respect to the standard structure on $\C^n \times
\C^n$). Using the implicit function theorem we study, 
for sufficiently small $\lambda$, the solutions of
$({\mathcal BP}_{\lambda})$ close to $(f^0,g^0)$. 
Similarly to Section~2.2 consider the map $u$
defined in a neighborhood of $(f^0,g^0,0)$ in 
$(\mathcal C^\alpha(\bar{\Delta}))^{4n} \times \mathbb R$ by:  

$$u: (f,g,\lambda) \mapsto \left(
\begin{matrix}
\zeta \in \partial \Delta \mapsto r(f(\zeta),\zeta^{-1}g(\zeta),\lambda)\\
\bar \partial f + q_1(\lambda,f)\partial f\\
\bar \partial g + q_2(\lambda,f) \partial g + q_3(\lambda,f)  g
\end{matrix}
\right).
$$
If $X:=(\mathcal C^ \alpha(\bar{\Delta}))^{4n}$ 
then its tangent map at $(f^0,g^0,0)$ has the form 

$$\begin{array}{llllll}
D_Xu(f^0,g^0,0): h=(h_1,h_2) \mapsto  \left(
\begin{matrix}
\zeta \in \partial \Delta \mapsto 
2 Re [G(f^0(\zeta),\zeta^{-1}g^0(\zeta),0)h]  \\
\bar \partial h_1 \\
\bar \partial h_2  
\end{matrix}
\right)
\end{array}$$
where for $\zeta \in \partial \Delta$ one has
$$\begin{array}{lllllll}
G(\zeta) = \left(
\begin{matrix}
\frac{\partial r_1}{\partial w_1}(f^0(\zeta),\zeta^{-1}g^0(\zeta),0)  &\cdots&\frac{\partial r_1}{\partial w_N}(f^0(\zeta),\zeta^{-1}g^0(\zeta),0)\\
\cdots&\cdots&\cdots\\
\frac{\partial r_N}{\partial w_1}(f^0(\zeta),\zeta^{-1}g^0(\zeta),0)&
\cdots&\frac{\partial r_N}
{\partial w_N}(f^0(\zeta),\zeta^{-1}g^0(\zeta),0)
\end{matrix}
\right)
\end{array}$$
with $N = 4n$ and $w = (z,t)$.

If the tangent map $D_Xu(f^0,g^0,0): (\mathcal A^{\alpha})^N
\longrightarrow (\mathcal 
C^{\alpha-1}(\Delta))^N$ is surjective  and has a
finite-dimensional kernel,  we may apply the implicit function theorem
as in Section~2.2 (see Proposition~\ref{tthh}) 
and conclude to the existence of a
finite-dimensional variety of nearby discs. 
In particular, consider the fibration $E$ over the disc $(f^0,g^0)$
with fibers
$E(\zeta) = \{ (z,t) : r^j(z,t,\zeta) = 0 \}$. Suppose that this
fibration is totally real. Then we have: 

\begin{proposition}\label{tthhh}
Suppose that the fibration $E$ is totally real. 
If the tangent map $D_Xu(f^0,g^0,0): (\mathcal A^{\alpha})^{4n}
\longrightarrow (\mathcal 
C^{\alpha-1}(\Delta))^{4n}$ is surjective 
 and has a finite-dimensional kernel,
then for every sufficiently small $\lambda$
the solutions of the boundary problem $({\mathcal BP}_{\lambda})$ form 
a smooth submanifold in the space $(\CC^{\alpha}(\Delta))^{4n}$.
\end{proposition}

In the next Section we present a sufficient condition for the 
surjectivity of the map $D_Xu(f^0,g^0,0)$. This is due to 
J.Globevnik~\cite{gl94,gl96} for the integrable case and
relies on the partial indices of the totally real fibration along $(f^0,g^0)$.

\section{Generation of stationary discs}

Let $D$ be a smoothly bounded domain in $\C^n$ with the
boundary $\Gamma$. According to \cite{le81} 
a continuous map $f:\bar\Delta \backslash \{ 0 \}
\rightarrow \bar D$, holomorphic on $\Delta \backslash \{ 0 \}$,
is called a {\it stationary} disc for $D$ (or
for $\Gamma$) if there
exists a holomorphic map
$\hat{f}: \Delta \backslash \{ 0 \} \rightarrow
T^*_{(1,0)}(\C^n)$, $\hat{f} \neq 0$, continuous on $\bar
\Delta \backslash \{ 0 \}$ and such that
\begin{itemize}
\item[(i)] $\pi \circ \hat{f} = f$
\item[(ii)] $\zeta \mapsto \zeta \hat{f}(\zeta)$ is in ${\mathcal O}(\Delta)$
\item[(iii)] $\hat{f}(\zeta) \in \Sigma_{f(\zeta)}(\Gamma)$ for every $\zeta$
in $\partial \Delta$.
\end{itemize}

We call $\hat{f}$ a {\it lift} of $f$ to the conormal bundle of $\Gamma$
(this is a meromorphic map from $\Delta$  into
$T^{*}_{(1,0)}(\C^n)$ whose values on the unit circle lie on 
$\Sigma(\Gamma)$).

We point out that originally Lempert gave this definition in a different form,
using the natural coordinates on the cotangent bundle of $\C^n$. The
present more geometric version in terms of the conormal
bundle is due to Tumanov~\cite{tu01}. This form is particularly useful
for our goals
since it can be transferred to the almost complex case.
\vskip 0,1cm
Let $f$ be a stationary disc for $\Gamma$.
It follows from Lemma~\ref{lemma} that if $\Gamma$ is a Levi nondegenerate
hypersurface, the conormal bundle $\Sigma(\Gamma)$ is a totally real
fibration along $f^*$. Conditions $(i)$, $(ii)$, $(iii)$ may be viewed as
a nonlinear boundary problem considered in Section~2.
If the associated tangent map is surjective, Proposition~\ref{tthhh} 
gives a description of all stationary
discs $\tilde{f}$ close to $f$, for a small deformation of $\Gamma$. 
When dealing with the standard complex structure on $\C^n$, the
bundle $T^*_{(1,0)}(\C^n)$  is a holomorphic vector bundle which can be 
identified, after projectivization of the fibers, with the
holomorphic bundle of complex hyperplanes that is with $\C^n \times
\mathbb P^{n-1}$. The conormal bundle $\Sigma(\Gamma)$ of a real
hypersurface $\Gamma$ may be naturally
identified, after this projectivization, with the bundle 
of holomorphic tangent spaces $H(\Gamma)$
over $\Gamma$. According to S.Webster
\cite{we78} this is a totally real submanifold in $\C^n \times \mathbb
P^{n-1}$. When dealing with the standard structure, we may therefore
work with projectivizations of lifts of stationary discs attached to
the holomorphic tangent bundle $H(\Gamma)$. The technical avantage
is that after such a projectivization lifts of stationary discs become
holomorphic, since the lifts have at most one pole of order 1 at the origin.
This idea was first used
by L.Lempert and then applied by several authors \cite{ba-le98,
ce95, sp-tr02}. 

When we consider almost complex 
deformations of the standard structure (and not just deformations of
$\Gamma$)
the situation is more complicated. The main
geometric difficulty is to prolong an almost complex structure $J$
from $\R^{2n}$ to the cotangent bundle $T^*(\R^{2n})$ in a certain
``natural'' way. As we will see later, such a prolongation is {\it
not} unique in contrast with the case of the integrable
structure. Moreover, if the cotangent bundle $T^*(\R^{2n})$ is
equipped with an almost complex structure, there is no natural
possibility to transfer this structure to the space obtained by the
projectivization of the fibers. Consequently we do not work with
projectivization of the cotangent bundle but we will deal with meromorphic
lifts of stationary discs. Representing such lifts $(\hat{f},\hat{g})$
in the form 
$(\hat{f},\hat{g})=(f,\zeta^{-1}g)$,
we will consider $(J_0,T_\lambda)$-holomorphic discs close to the 
$(J_0,T_0)$-holomorphic disc $(f,g)$. The disc $(f,g)$
satisfies a nonlinear boundary problem of 
Riemann-Hilbet type $(\mathcal BP_\lambda)$.
When an almost complex structure on the
cotangent bundle is fixed, we may view it as an elliptic prolongation
of an initial almost complex structure on the base and apply the
implicit function theorem as in previous section. This avoids
difficulties coming from the projectivization of almost complex fibre spaces. 

\subsection{Maslov index and Globevnik's condition}
We denote by $GL(N,\C)$
the group of invertible $(N \times N)$ complex matrices
and by $GL(N,\R)$ the group of all such matrices with real entries.
Let $0 < \alpha < 1$ and let
$B:\partial \Delta \rightarrow GL(N,\C)$ be of class $\mathcal C^\alpha$. 
According to \cite{ve67} (see also \cite{cl-go81}) $B$ admits the factorization
$B(\tau) = F^{+}(\tau)\Lambda(\tau)F^{-}(\tau), \tau \in \partial \Delta$,
where:

$\bullet$ $\Lambda$ is a diagonal matrix of the form
$\Lambda(\tau) = diag(\tau^{k_1},\dots,\tau^{k_N})$, 

$\bullet$ $F^{+}: \bar{\Delta} \rightarrow GL(N,\C)$ is
of class $\mathcal C^\alpha$ on $\bar{\Delta}$ and holomorphic in $\Delta$,

$\bullet$ $F^{-}: [\C \cup \{ \infty \}] \backslash \Delta
\rightarrow GL(N,\C)$
is of class $\mathcal C^\alpha$ on $[\C \cup \{ \infty \}] \backslash \Delta$
and holomorphic on
$[\C \cup \{ \infty \}] \backslash \bar \Delta$.

\vskip 0,1cm The integers $k_1 \geq \cdots \geq k_n$ are called the partial 
indices of $B$. 

\vskip 0,1cm
Let $E$ be a totally real fibration over the unit circle. For every 
$\zeta \in \partial \Delta$ consider the ``normal''
vectors $\nu_j(\zeta) =
(r^j_{\bar{z}^1}(f(\zeta),\zeta),\dots,
r^j_{\bar{z}^N}(f(\zeta),\zeta))$, $j=1,\dots,N$.
We denote by $K(\zeta) \in GL(N,\C)$ the matrix with rows
$\nu_1(\zeta),\dots,\nu_N(\zeta)$ and we set
 $B(\zeta) := -\overline{K(\zeta)}^{-1}K(\zeta)$,
$\zeta \in \partial \Delta$. The partial indices of the map
$B: \partial \Delta \rightarrow GL(N,\C)$
are called {\it the partial indices} of the fibration ${E}$
along the disc $f$ and their sum is
called the {\it total index} or the {\it Maslov index of ${ E}$ along $f$}. 
The following result is due to
J. Globevnik~\cite{gl96}:

\vskip 0,1cm
\noindent{\bf Theorem} : {\it Suppose that all the partial indices of
the totally real fibration ${E}$ along $f$ are $\geq -1$ 
and denote by $k$ the Maslov index of $E$ along 
$f$. Then the linear map from $(\mathcal A^{\alpha})^N$ to $(\mathcal 
C^{\alpha-1}(\Delta))^N$ given by 
$h \mapsto 2 Re [G h]$
is surjective and has a $k$ dimensional kernel.}

\vskip 0,1cm
Proposition~\ref{tthh} may be restated in terms of partial indices
as follows~:

\begin{proposition}\label{tthh1}
Let $f^0:\bar \Delta \rightarrow \C^N$ be a $J_0$-holomorphic
disc attached to a totally real fibration ${E}$ in $\C^N$. Suppose
that all the partial indices of ${E}$ along $f^0$ are $\geq
-1$. Denote by $k$ the Maslov index of ${E}$ along $f^0$. Let also
$E_\varepsilon$ be a smooth totally real deformation of ${E}$ and
$J_\lambda$ be a smooth almost complex deformation of $J_0$ in a
neighborhood of $f(\bar{\Delta})$. Then there exists $\delta_0,
\varepsilon_0, \lambda_0 >0$ such that for every $0 \leq \varepsilon
\leq \varepsilon_0$ and for every $0 \leq \lambda \leq \lambda_0$  
the set of $J_\lambda$-holomorphic discs attached to $E_\varepsilon$
and such that $\parallel f -f^0 \parallel_{\alpha} \leq \delta_0$ forms
a smooth $(N+k)$-dimensional submanifold 
$\mathcal A_{\varepsilon,\lambda}$ in the Banach space 
$(\mathcal C^\alpha(\bar{\Delta}))^N$.
\end{proposition}

Globevnik's result was applied to the study of stationary discs in
some classes of domains in $\C^n$ by M.Cerne~\cite{ce95} and
A.Spiro-S.Trapani~\cite{sp-tr02}.  Since they worked with the
projectivization of the conormal bundle, we explicitely compute, for
reader's convenience and completeness of exposition, partial indices
of {\it meromorphic} lifts of stationary discs for the unit sphere in
$\C^n$.

\vskip 0,1cm
Consider the unit sphere $\Gamma:=\{z \in \C^n : 
z^1\bar{z}^1 + \cdots +z^{n}\bar{z}^{n} -1 = 0\}$ in $\C^n$. 
The conormal bundle $\Sigma(\Gamma)$ is given in the $(z,t)$ coordinates 
 by the equations 
$$
(S)\left\{
\begin{array}{cll}
z^1\bar{z}^1 + \cdots +z^{n}\bar{z}^{n} - 1 = 0, & &\\
t_1 = c \bar{z}^1,\dots,t_{n} = c \bar{z}^{n}, &c \in \mathbb R.&
\end{array}
\right.
$$
According to \cite{le81}, every stationary disc for $\Gamma$
 is extremal for the Kobayashi metric. Therefore,
such a stationary disc $f^0$ centered at the origin is linear by the
Schwarz lemma. So, up to a unitary transformation, we have
 $f^0(\zeta) = (\zeta,0,\dots,0)$ with lift 
$(\widehat{f^0},\widehat{g^0})(\zeta) =
(\zeta,0,\dots,0,\zeta^{-1},0,\dots,0)=(f^0,\zeta^{-1}g^0)$ to the 
conormal bundle.
Representing nearby meromorphic discs in the form
$(z,\zeta^{-1}w)$ and eliminating the parameter $c$ in system 
$(S)$ we obtain that holomorphic discs
$(z,w)$ close to $(f^0,g^0)$ satisfy for $\zeta \in \partial \Delta$:

$$
\begin{array}{lll}
r^1(z,w) & = & z^1\bar{z}^1 + \cdots +z^{n}\bar{z}^{n} - 1 = 0,\\
r^2(z,w) & = & i z^1 w_1\zeta^{-1} - i\bar{z}^1\bar{w}_1\zeta = 0,\\
r^3(z,w) & = & \bar{z}^1w_2\zeta^{-1} - \bar{z}^2w_1\zeta^{-1} + 
z^1\bar{w}_2 \zeta-
z^2\bar{w}_1\zeta = 0,\\
r^4(z,w) & = & i\bar{z}^1w_2\zeta^{-1} - i\bar{z}^2w_1\zeta^{-1} -
 iz^1\bar{w}_2\zeta +
iz^2\bar{w}_1\zeta = 0,\\
r^5(z,w) & = & \bar{z}^1w_3\zeta^{-1} - \bar{z}^3w_1\zeta^{-1} +
 z^1\bar{w}_3 \zeta-
z^3\bar{w}_1\zeta = 0,\\
r^6(z,w) & = & i\bar{z}^1w_3\zeta^{-1} - i\bar{z}^3w_1\zeta^{-1} - 
iz^1\bar{w}_3\zeta +
iz^3\bar{w}_1\zeta = 0,\\
& &\cdots\\
r^{2n-1}(z,w) & = & \bar{z}^1w_n\zeta^{-1} - \bar{z}^nw_1\zeta^{-1} + 
z^1\bar{w}_n \zeta-
z^n\bar{w}_1\zeta = 0,\\
r^{2n}(z,w) & = & i\bar{z}^1w_n\zeta^{-1} - i\bar{z}^nw_1\zeta^{-1} - 
iz^1\bar{w}_n\zeta +
iz^n\bar{w}_1\zeta = 0.
\end{array}
$$

Hence the $(2n \times 2n)$-matrix $K(\zeta)$ has the following expression:
$$
\left( 
\begin{matrix}
\zeta&0&0& \cdots &0&0&0&0& \cdots & 0\\
-i\zeta&0&0& \cdots &0&-i&0&0&\cdots &0\\
0&-\zeta^{-1}&0& \cdots &0&0&\zeta^2&0&\cdots &0\\
0&-i\zeta^{-1}&0& \cdots &0&0&-i\zeta^2&0&\cdots &0\\
0&0&-\zeta^{-1}& \cdots &0&0&0&\zeta^2&\cdots &0\\
0&0&-i\zeta^{-1}& \cdots &0&0&0&-i\zeta^2&\cdots &0\\
\cdots & \cdots & \cdots & \cdots & 
\cdots & \cdots & \cdots & \cdots & \cdots &\cdots\\
0&0&0& \cdots &-\zeta^{-1}&0&0&0& \cdots &\zeta^2\\
0&0&0& \cdots &-i\zeta^{-1}&0&0&0& \cdots &-i\zeta^2
\end{matrix}
\right)
$$
and a direct computation shows that $-B =  \bar{K}^{-1} K$ has
the form 
$$
\left( \begin{matrix}C_1&C_2\\
C_3&C_4
\end{matrix}
\right),
$$
where the $(n \times n)$ matrices $C_1,\dots,C_4$ are given by

$$
C_1 = \left( \begin{matrix}\zeta^2&0&.&0\\
0&0&.&0\\
0&0&.&0\\
.&.&.&.\\
0&0&.&0
\end{matrix}
\right), \ \ C_2 = \left( \begin{matrix}0&0&.&0\\
0&-\zeta&.&0\\
0&0&-\zeta&0\\
.&.&.&.\\
0&0&0&-\zeta
\end{matrix}
\right),
$$
$$
C_3 = \left( \begin{matrix}-2\zeta&0&.&0\\
0&-\zeta&.&0\\
0&0&.&0\\
.&.&.&.\\
0&0&.&-\zeta
\end{matrix}
\right), \ \ C_4 = \left( \begin{matrix}-1&0&.&0\\
0&0&.&0\\
0&0&.&0\\
.&.&.&.\\
0&0&.&0
\end{matrix}
\right).
$$
\vskip 0,1cm

We point out that the matrix 

$$\left( \begin{matrix}\zeta^2&0\\
\zeta&1
\end{matrix}
\right)$$
admits the following factorization:

$$\left( \begin{matrix}1&\zeta\\
0&1
\end{matrix} \right ) \times \left( \begin{matrix}-\zeta&0\\
0&\zeta
\end{matrix} 
\right) \times 
\left( \begin{matrix}0&1\\
1&\zeta^{-1}
\end{matrix} \right).$$

Permutating the lines (that is multiplying $B$ by some nondegenerate matrics 
with constant coefficients) and using the above factorization of $(2 \times
2)$ matrices, we obtain the following 

\begin{proposition}\label{PPRRR}
All the partial indices of the conormal bundle of the unit sphere along a 
meromorphic disc of a stationary disc are equal to one 
and the Maslov index is equal to $2n$.
\end{proposition}

Proposition~\ref{PPRRR} enables to apply Proposition~\ref{tthh}
to construct the family of stationary discs attached to the unit
sphere after a small deformation of the complex structure. Indeed
denote by $r^j(z,w,\zeta,\lambda)$ $\CC^{\infty}$-smooth functions 
coinciding for
$\lambda = 0$ with the above functions $r^1,\dots,r^{2n}$.  

\vskip 0,1cm
In the end of this Subsection we
make the two following assumptions:
\begin{itemize}
\item[(i)] $r^1(z,w,\zeta,\lambda) =  z^1\bar{z}^1 +
\cdots+z^{n}\bar{z}^{n} - 1$, meaning that the sphere is not deformed
\item[(ii)] $r^j(z,tw,\zeta,\lambda)
= tr^j(z,w,\zeta,\lambda)$ for every $j \geq 2, \ t \in \R$.
\end{itemize}

Geometrically this means that given $\lambda$, the set $\{ (z,w):
r^j(z,w,\lambda) = 0 \}$ is a real vector bundle with one-dimensional
fibers over the unit sphere.

Consider an almost
complex deformation $J_{\lambda}$ of the standard structure on $\mathbb B_n$
and an
elliptic prolongation $T_{\lambda}$ of $J_{\lambda}$ on $\mathbb B_n \times 
\mathbb R^{2n}$. Recall that the
corresponding $\bar\partial_{T_{\lambda}}$-equation is given by the 
following system
$$
({\mathcal E})
\left\{
\begin{array}{lll}
& &\bar \partial f + q_1(\lambda,f))\partial f = 0,\\
& &\bar \partial g + q_2(\lambda,f))\partial g + q_3(\lambda,f) g = 0.
\end{array}
\right.
$$
Consider now the corresponding boundary problem:
$$
({\mathcal BP}_{\lambda})
\left\{
\begin{array}{lll}
& &r(f,g,\zeta,\lambda) = 0, \zeta \in \partial \Delta,\\
& &\bar \partial f + q_1(\lambda,f)\partial f = 0,\\
& &\bar \partial  g + q_2(\lambda,f)\partial g +
q_3(\lambda,f) g = 0.
\end{array}
\right.
$$

Combining Proposition~\ref{PPRRR} with the results of Section~2,
we obtain the following
\begin{proposition}\label{THEO}
For every sufficiently small positive $\lambda$, the set of solutions of
 $({\mathcal BP}_{\lambda})$, close enough to the disc
$(\widehat{f^0},\widehat{g^0})$, forms a smooth $4n$-dimensional
submanifold $V_{\lambda}$ in the space
$\CC^{\alpha}(\bar \Delta)$ (for every noninteger $\alpha > 1$).
\end{proposition}
Moreover, in view of the assumption (ii) and of
the linearity of $({\mathcal E})$ with respect to the fiber component $g$,
we also have the following

\begin{corollary}\label{CORO}
The projections of discs from
$V_{\lambda}$
to the base $(\R^{2n},J_{\lambda})$ form a $(4n-1)$-dimensional subvariety in
$\CC^{\alpha}(\bar\Delta)$.
\end{corollary}

Geometrically the solutions $(f,g)$ of the boundary problem
$({\mathcal BP}_{\lambda})$ are
such that the discs $(f,\zeta^{-1}g)$ are attached to the conormal
bundle of the unit sphere with respect to the standard structure. In
particular, if $\lambda = 0$ then every such disc satisfying $f(0) =
0$ is linear.

\section{Foliation and ``Riemann map'' associated with an elliptic
prolongation of an almost complex structure}
In this Section we study the geometry of stationary discs in the unit ball 
after
a small almost complex perturbation of the standard structure. 
The idea is simple since these discs are small deformations of the 
complex lines passing through the origin in the unit ball.  
\subsection{Foliation associated with an elliptic prolongation} 
Fix a vector $v^0$ with $\vert\vert v^0 \vert\vert = 1$ and consider the
corresponding stationary disc $f^0:\zeta \mapsto \zeta v^0$. Denote by
$(\widehat{f^0},\widehat{g^0})$ its lift to the conormal bundle of
the unit sphere. 
Consider a smooth deformation $J_{\lambda}$ of the standard structure
on the unit ball ${\mathbb B_n}$ in $\C^n$.
For sufficiently small $\lambda_0$ fix a prolongation $T_{\lambda}$ of the
structure $J_{\lambda}$ on $\mathbb B_n \times \mathbb R^{2n}$, 
where $\lambda \leq
\lambda_0$. Then the solutions of the associated boundary
problem $({\mathcal BP}_{\lambda})$  form a $4n$-parameter family
of $(J_0,T_\lambda)$-holomorphic maps from $\Delta \backslash \{0\}$ to
$\C^{n} \times \C^{n}$. Given such a solution $(f^\lambda,g^\lambda)$, 
consider the disc
$(\widehat{f^\lambda},\widehat{g^\lambda}) := 
(f^\lambda, \zeta^{-1}g^\lambda)$. 
In the case where $\lambda = 0$ 
this is just the lift of a stationary disc for the unit
sphere to its conormal bundle. The set of solutions of the problem
$(\mathcal BP_\lambda)$ (considered in Section~3.1), 
close to $(\widehat{f^0},\widehat{g^0})$, forms a
smooth submanifold of real dimension $4n$ 
in $(\mathcal C^\alpha(\bar{\Delta}))^{4n}$ according
to Theorem~\ref{THEO}. Hence there is a neighborhood $V_0$ of $v^0$ in
$\mathbb R^{2n}$ and a smooth real hypersurface $I_{v^0}^\lambda$ in $V_0$ 
such that for every $\lambda \leq \lambda_0$ and for every
$v \in I_{v^0}^\lambda$ there is one and only one solution 
$(f^\lambda_v,g^\lambda_v)$ of $(\mathcal BP_\lambda)$, up to multiplication
of the fiber component $g^\lambda_v$ by a real constant,
such that 
$f^\lambda_v(0) = 0$ and $df^\lambda_v(0)(\partial / \partial x) = v$.

We may therefore consider the map
$$
F_0^{\lambda}: (v,\zeta) \in I_{v^0}^\lambda \times \bar{\Delta} \mapsto
(f^\lambda_v,g^\lambda_v)(\zeta).
$$

This is a smooth map with respect to $\lambda$ 
close to the origin in $\mathbb R$.

Denote by $\pi$ the canonical
projection $\pi: \mathbb B_n \times \R^{2n} \rightarrow \mathbb B_n$
and consider the composition $\widehat{F_0^{\lambda}} = \pi \circ
F_0^{\lambda}$. 
This is a smooth map defined for $0 \leq \lambda <
\lambda_0$
and such that

\begin{itemize}
\item[(i)] $\widehat{F_0^{0}}(v,\zeta) 
= v\zeta$, for every $\zeta \in \bar\Delta$ and for every $v \in I_{v^0}^\lambda$.
\item[(ii)] For every $\lambda \leq \lambda_0$, 
$\widehat{F_0^\lambda}(v,0) = 0$.
\item[(iii)] For every fixed $\lambda \leq \lambda_0$ and every 
$v \in I_{v^0}^\lambda$ the map 
$\widehat{F_0^{\lambda}}(v,\cdot)$ is
a $(J_0,J_{\lambda})$-holomorphic disc attached to the unit sphere.
\item[(iv)] For every fixed $\lambda$, different values of
$v \in I_{v^0}^\lambda$ define different discs.
\end{itemize}

\begin{definition} We call the family 
$(\widehat{F_0^\lambda}(v,\cdot))_{v \in I_{v^0}^\lambda}$ {\it canonical }
discs associated with the boundary problem $({\mathcal
BP}_{\lambda})$.
\end{definition}

 We stress that by 
a canonical disc {\it we always mean a disc centered at the
origin}. The preceding condition $(iv)$ may be restated as follows:

\begin{lemma}\label{LEMMA1} 
For $\lambda < \lambda_0$ every canonical disc is uniquely 
determined by its tangent vector at the origin.
\end{lemma}

In the next Subsection we glue the sets $I_{v}^\lambda$, depending on vectors
$v \in \mathbb S^{2n-1}$, to define the global indicatrix of $F^\lambda$.

\subsection{Indicatrix} For $\lambda < \lambda_0$ consider canonical discs
in $(\mathbb B_n,J_{\lambda})$ centered at the origin and admitting lifts
close to $(\widehat{f^0},\widehat{g^0})$. 

As above we denote by $I_{v^0}^{\lambda}$ 
the set of tangent vectors at the origin of
canonical discs whose lift is close to $(\widehat{f^0},\widehat{g^0})$. 
Since these vectors depend smoothly on
parameters $v$ close to $v^0$ and $\lambda \leq \lambda_0$, 
$I_{v^0}^{\lambda}$ is a smooth deformation of a piece of
the unit sphere $\mathbb S^{2n-1}$. So this is a smooth real
hypersurface in $\C^n$ in a neighborhood of $v^0$.  
Repeating this construction for every vector $v \in \mathbb S^{2n-1}$ 
we may find a finite covering of $\mathbb S^{2n-1}$ by open
connected sets $U_j$ such that for every $j$ the nearby stationary
discs with tangent vectors at the origin close to $v$ are given by
$\widehat{F^{\lambda}_j}$. Since every nearby stationary disc is uniquely
determined by its tangent vector at the origin, we may glue the maps
$\widehat{F^{\lambda}_j}$ to the map $\widehat{F^{\lambda}}$ 
defined for every $v
\in \mathbb S^{2n-1}$ and every $\zeta \in \bar{\Delta}$. The tangent vectors
of the constructed family of stationary discs form a smooth real
hypersurface $I^{\lambda}$ which is a small deformation of the unit
sphere. This hypersurface is an analog of the indicatrix for the
Kobayashi metric (more precisely, its boundary).

We point out that the local indicatrix $I_{v^0}^\lambda$ 
for some fixed $v^0 \in \mathbb S^{2n-1}$ is also useful (see Section 5).

\subsection{Circled property and Riemann map}
If $\lambda$  is  small enough, the hypersurface 
$I^{\lambda}$ is strictly pseudoconvex with respect to the standard
structure. Another important property of the ``indicatrix'' is its
invariance with respect to the linear action of the unit circle.

Let $\lambda \leq \lambda_0$, $v \in I^\lambda$ and 
$f_v^\lambda:=\widehat{F^\lambda}(v,\cdot)$. For $\theta \in \mathbb R$ 
we denote by
$f_{v,\theta}^\lambda$ the $J_\lambda$-holomorphic disc in $\mathbb B_n$
defined by $f^\lambda_{v,\theta} : \zeta \in \Delta \mapsto
f_v^\lambda(e^{i\theta}\zeta)$. We have~:

\begin{lemma}\label{LEMMA2}
For every $0 \leq \lambda < \lambda_0$, every $v \in I^\lambda$
and every $\theta \in \R$ we have~: $f_{v,\theta}^\lambda \equiv
f_{e^{i\theta}v}^\lambda$.
\end{lemma}

\proof Since $f_v^\lambda$ is a canonical disc, the disc
$f_{v,\theta}^\lambda$ has a  lift   close to the lift
of the disc $\zeta \mapsto e^{i\theta}v\zeta$. Then according to Lemma
\ref{circle} 
$f_{v,\theta}^\lambda$ is a canonical disc close to the
linear disc $\zeta \mapsto e^{i\theta}v\zeta$. Since the first jet
of $f_{v,\theta}^\lambda$ coincides with the first jet of
$f_{e^{i\theta}v}^\lambda$, these two nearby stationary discs coincide
according to Lemma~\ref{LEMMA1}. \qed

This statement implies that for any $w \in I^{\lambda}$ the vector
$e^{i\theta}w$ is in $I^{\lambda}$ as well. 

It follows from the
above arguments that there exists a natural parametrization of the set
of canonical discs by their tangent vectors at the origin,
that is by the points of $I^{\lambda}$. The map
$$
\begin{array}{llcll}
\widehat{F^\lambda} &:&  I^\lambda \times \Delta 
& \rightarrow & \mathbb B_n\\
    
    & &  (v,\zeta) & \mapsto     & f_v^\lambda(\zeta)
\end{array}
$$ 
is smooth on $I^\lambda \times \Delta$. Moreover, if we fix a
small positive constant $\varepsilon_0$, then by shrinking $\lambda_0$
if necessary there is, for every $\lambda < \lambda_0$, a smooth
function $\widehat{G^\lambda}$ defined on $I^\lambda \times \Delta$,
satisfying $\|\widehat{G^\lambda(v,\zeta)}\| \leq \varepsilon_0 |\zeta|^2$ on
$I^\lambda \times \Delta$, such that for every $\lambda <
\lambda_0$ we have on $I^\lambda \times \Delta$~:

\begin{equation}\label{equation}
\widehat{F^ \lambda}(v,\zeta) = \zeta v + \widehat{G^\lambda}(v,\zeta).
\end{equation}

\vskip 0,1cm
Consider now the restriction of $\widehat{F^\lambda}$ to $I^\lambda
\times [0,1]$. This is a smooth map, still denoted by 
$\widehat{F^\lambda}$.  We
have the following~:
\begin{proposition}\label{PROPRO}
There exists $\lambda_1 \leq \lambda_0$ such that 
for every $\lambda < \lambda_1$ the family
$(\widehat{F^\lambda}(v,r))_{(v,r) \in I^\lambda \times [0,1[}$
is a real foliation of $\mathbb B_n \backslash \{0\}$. 
\end{proposition}
\proof 

\noindent{\bf Step 1}.
For $r \neq 0$ we write $w:=rv$. Then $r=\|w\|$, $v=w/\|w\|$ 
and we denote by $\widetilde{F^\lambda}$ the function 
$\widetilde{F^\lambda}(w):=\widehat{F^\lambda}(v,r)$. 
For $\lambda < \lambda_0$, $\widetilde{F^\lambda}$ is a smooth map 
of the variable $w$ on
$\mathbb B_n \backslash\{0\}$, satisfying~:
$$
\widetilde{F^\lambda}(w) = w + \widetilde{G^\lambda}(w)
$$
where $\widetilde{G^\lambda}$ is a smooth map on $\mathbb B_n \backslash \{0\}$
with $\|\tilde{G}^\lambda(w)\| \leq \varepsilon_0\|w\|^2$ on $\mathbb B_n
\backslash\{0\}$. This implies that $\widetilde{F^\lambda}$ is a local
diffeomorphism at each point in $\mathbb B_n \backslash\{0\}$, and so that
$\widehat{F^\lambda}$ is a local diffeomorphism at each point in
$I^\lambda \times ]0,1[$. Moreover, the
condition $\|\widetilde{G^\lambda}(w)\| \leq \varepsilon_1\|w\|^2$ on
$\mathbb B_n \backslash\{0\}$ implies that $\widetilde{G^\lambda}$ is
differentiable at the origin with $d\widetilde{G^\lambda}(0) = 0$. Hence
by the implicit function theorem there exists $\lambda_1 < \lambda_0$
such that the map $\widetilde{F^\lambda}$ is a local diffeomorphism at the
origin for $\lambda < \lambda_1$. So there exists $0<r_1<1$ and a 
neighborhood $U$ of the origin in $\C^n$ such that
$\widehat{F^\lambda}$ is a diffeomorphism from $I^\lambda \times
]0,r_1[$ to $U \backslash \{0\}$, for $\lambda < \lambda_1$.

\vskip 0,1cm
\noindent{\bf Step 2}. We show that $\widehat{F^\lambda}$ 
is injective on $I^\lambda \times ]0,1]$ for sufficiently small $\lambda$. 
Assume by contradiction that for every $n$ there exist 
$\lambda_n \in \mathbb R$, $r_n, r'_n \in ]0,1]$, $v^n, w^n \in I^{\lambda_n}$ 
such that:

$ \bullet \lim_{n \rightarrow \infty}\lambda_n =0, 
\ \lim_{n \rightarrow \infty}r_n=r,\ 
\lim_{n \rightarrow \infty}r'_n=r',
$

$\bullet
\lim_{n \rightarrow \infty}v^n=v \in \mathbb S^{2n-1}, \ 
\lim_{n \rightarrow \infty}w^n=w \in \mathbb S^{2n-1}
$

and satisfying
$$
\widehat{F^{\lambda_n}}(v^n,r_n) = \widehat{F^{\lambda_n}}(w^n,r'_n)
$$
for every $n$. Since $\widehat{F}$ is smooth with respect to $\lambda, v, r$,
it follows that $\widehat{F^0}(v,r) = \widehat{F^0}(w,r')$ and so
$v=w$ and $r=r'$. If $r < r_1$ then the contradiction follows from the fact 
that $\widehat{F^\lambda}$ is a diffeomorphism from 
$I^\lambda \times ]0,r_1[$ to $U \backslash \{0\}$. If $r \geq r_1$
then for every neighborhood $U_\infty$ of $rv$ in 
$\mathbb B_n \backslash \{0\}$, $r_nv^n \in U_\infty$ and 
$r_n'w^n \in U_\infty$ for sufficiently large $n$. Since we may choose 
$U_\infty$ such that $\widehat{F^\lambda}$ is a diffeomorphism from a
neighborhood of $(v,r)$ in $I^\lambda \times ]r_1,1]$ uniformly with respect
to $\lambda <<1$, we still obtain a contradiction.

\vskip 0,1cm
\noindent{\bf Step 3}. We show that $\widehat{F^\lambda}$ is surjective
from $I^\lambda \times ]0,1[$ to $\mathbb B_n \backslash \{0\}$. 
It is sufficient to show that $\widehat{F^\lambda}$ is surjective from
$I^\lambda \times [r_1,1[$ to $\mathbb B_n \backslash U$.
Consider the nonempty set 
$E_\lambda=\{w \in \mathbb B_n \backslash U : w=\widehat{F^\lambda}(v,r) \ 
{\rm for \ some \ }(v,r) \in I^\lambda \times ]r_1,1[\}$. Since the jacobian
of $\widehat{F^\lambda}$ does not vanish for $\lambda=0$ and 
$\widehat{F^\lambda}$ is smooth with respect to $\lambda$, 
the set $E_\lambda$ is open for
sufficiently small $\lambda$.
Moreover it follows immediately from its definition that $E_\lambda$ 
is also closed in $\mathbb B_n \backslash U$. Thus 
$E_\lambda=\mathbb B_n \backslash U$.

These three steps prove the result. \qed

\vskip 0,1cm
We can construct now the map $\Psi_{T_\lambda}$ for $\lambda <
\lambda_1$. For every $z \in \mathbb B_n \backslash \{0\}$ consider the unique
couple $(v(z),r(z)) \in I^\lambda \times ]0,1[$ such that
$f_v^\lambda$ is the unique canonical disc passing through $z$ (its
existence and unicity are given by Proposition~\ref{PROPRO}) with
$f_{v(z)}^\lambda(0) = 0$, $df_{v(z)}^\lambda(0)(\partial / \partial
x) = v(z)$ and $f_{v(z)}^\lambda(r(z)) = z$. The map $\Psi_{T_\lambda}$ is
defined by~:
$$
\begin{array}{llcll}
\Psi_{T_\lambda} &: & \bar{\mathbb B}_n \backslash \{0\} & \rightarrow & \C^n\\
     &  &  z                   & \mapsto     & r(z) v(z).
\end{array}
$$
\begin{definition}
The map $\Psi_{T_\lambda}$ is called the Riemann map associated with a
prolongation $T_{\lambda}$ of the almost complex structure
$J_{\lambda}$.
\end{definition}
This map is an analogue of the circular representation of a strictly
convex domain introduced by L.Lempert~\cite{le81}. The term ``Riemann map''
was used by S. Semmes~\cite{se92} for a slightly different map where the
vector $v(z)$ is normalized (and so such a map takes values in the unit ball).
In this paper we work with the indicatrix since this is more convenient for
our applications.

The Riemann map $\Psi_{T_\lambda}$ has the following properties~:
\begin{proposition}\label{PROPRO2}
\noindent $(i)$ For every $(v,\zeta) \in I^\lambda \times
\Delta$ we have $(\Psi_{T_\lambda} \circ f_v^\lambda)(\zeta) =
\zeta v$ and so $\log\|(\Psi_{T_\lambda} \circ f_v^\lambda)(\zeta)\| =
\log|\zeta|$.

\noindent $(ii)$ There exist constants $0 < C' < C$ such that $C'
\|z\| \leq \|\Psi_{T_\lambda}(z)\| \leq C \|z\|$ on $\mathbb B_n$.
\end{proposition}

\vskip 0,1cm
\proof $(i)$ Let $\zeta = e^{i\theta}r \in \Delta(0,r_0)$ with $\theta
\in [0,2\pi[$. Then $f_v^\lambda(\zeta) = f_v^\lambda(e^{i\theta}r) =
f_{e^{i\theta}v}^\lambda(r)$. Hence we have $(\Psi_{T_\lambda} \circ
f_v^\lambda)(\zeta) = \Psi_{T_\lambda}( f_{e^{i\theta}v}^\lambda(r)) =
e^{i\theta}vr = \zeta v$.

$(ii)$ Let $z \in \mathbb B_n \backslash \{0\}$. Then according to
equation~(\ref{equation}) we have the inequality
$\|\Phi^{\lambda}(z)\|\ (1-\varepsilon_1\|\Psi_{T_\lambda}(z)\|) \leq
\|z\|^2 \leq \Psi_{T_\lambda}(z) \ (1+ \varepsilon_1
\|\Psi_{T_\lambda}(z)\|)$. Since $\|\Psi_{T_\lambda}(z)\| \leq 1$ we obtain
the desired inequality with $c'=1/1+\varepsilon_1$ and
$c=1/1-\varepsilon_1$. \qed

\vskip 0,1cm
From the above analysis it follows the basic properties of the Riemann map.

\begin{proposition}\label{PR1}

\begin{itemize}
\item[]
\item[(i)] The indicatrix $I^{\lambda}$ is a compact circled smooth
$J_{\lambda}$-strictly pseudoconvex hypersurface bounding a domain
denoted by $\Omega^{\lambda}$.
\item[(ii)] The Riemann map $\Psi_{T_\lambda}: \bar {\mathbb B}_n \backslash
\{ 0 \} \rightarrow \bar\Omega^{\lambda} \backslash \{ 0 \}$ is a smooth
diffeomorphism.
\item[(iii)] For every canonical disc $f_v^{\lambda}$ we have
$\Psi_{T_{\lambda}}\circ f^{\lambda}_v(\zeta) = v\zeta$.
\end{itemize}
\end{proposition}

We note that the Riemann map possesses further important structure
properties depending on the choice of a prolongation $T_{\lambda}$
of an almost complex structure $J_{\lambda}$.  

\subsection{Local Riemann map}
In Subsection~4.2 we introduced the notion of local indicatrix 
$I_{v^0}^\lambda$ for $v^0 \in \mathbb S^{2n-1}$. We may localize the
notion of the Riemann map, introducing a similar associated with the
local indicatrix. Denote by $\Omega_{v^0}^\lambda$ the set 
$I_{v^0}^\lambda \times
[0,1[$. The arguments used in the proof of Proposition~\ref{PROPRO}
show that $\widehat{F^\lambda}(\Omega_{v^0}^\lambda)$ 
is foliated by stationary discs centered at the origin.
We may therefore define the Riemann map $\Psi_{T_\lambda,v^0}$ on 
$\widehat{F^\lambda}(\Omega_{v^0}^\lambda)$ by:
$$
\Psi_{T_\lambda,v^0}(z)=r(z)v(z)
$$
where $v(z)$ is the tangent vector at the origin of the unique stationary disc
$f_{v(z)}^\lambda$ passing through $z$ and $f_{v(z)}^\lambda(r(z)) = z$.

\begin{remark}\label{REM}
We point out that the Riemann map can be defined in any sufficiently small 
deformation of the unit ball and satisfies all the same properties.
\end{remark}
\section{ Riemann map and local geometry of almost complex manifolds}
As we have seen previously, a choice of an elliptic prolongation
$T_{\lambda}$ of an almost complex stucture $J_{\lambda}$ on a vector
bundle over the ball (for
$\lambda$ small enough) allows to define the foliation by canonical
discs and
the Riemann map. Further properties of such a map depend on the choice
of $T_{\lambda}$ which is not unique. In this section we consider
a natural definition of $T_{\lambda}$ as the vertical lift of the almost
complex structure $J_\lambda$. This construction is well-known in
differential geometry on the tangent and cotangent bundles, see~\cite{is-ya73}.
For reader's convenience we recall it.

\subsection{Vertical lift to the cotangent bundle} Our goal
now is to construct a family of Lempert discs invariant with respect
to biholomorphic transformations of an almost complex manifold with
boundary. We recall the definition of the canonical lift of an almost
complex structure $J$ on $M$ to the cotangent bundle $T^*M$. Set $m=
2n$. We use the following notations. Suffixes A,B,C,D take the values
$1$ to $2m$, suffixes $a,b,c,\dots$,$h,i,j,\dots$ take the values $1$ to
$m$ and $\bar j = j+ m$, $\dots$ The summation notation for
repeated indices is used. If the notation $(\varepsilon_{AB})$,
$(\varepsilon^{AB})$, $(F_{B}^{{}A})$ is used for matrices, the suffix
on the left indicates the column and the suffix on the right indicates
the row. We denote local coordinates on $M$ by $(x^1,\dots,x^n)$ and by
$(p_1,\dots,p_n)$ the fiber coordinates.

Recall that the cotangent space $T^*(M)$ of $M$ possesses the {\it
canonical contact form} $ \theta$ given in local coordinates by

$$\theta = p_idx^i.$$
The cotangent lift $\varphi^*$ of any diffeomorphism $\varphi$ of $M$
is contact with respect to $\theta$, that is $\theta$ does not depend on
the choice of local coordinates on $T^*(M)$. 

The exterior derivative $d\theta$ of $\theta$ defines {\it the
canonical
symplectic structure} of $T^*(M)$:

$$d\theta = dp_i \wedge dx^i$$
which is also independent of local coordinates in view of the
invariance of the exterior derivative. Setting $d\theta =
(1/2)\varepsilon_{CB}dx^C \wedge dx^b$ (where $dx^{\bar j} =
dp_j$), we have

$$
(\varepsilon_{CB}) = \left(
\begin{matrix}
0 & I_n\\
-I_n & 0
\end{matrix}
\right).
$$

Denote by $(\varepsilon^{BA})$ the inverse matrix, we write
$\varepsilon^{-1}$ for the tensor field of type (2,0) whose component
are $(\varepsilon^{BA})$. By construction, this definition does not
depend on the choice of local coordinates.

Let now $E$ be a tensor field of type (1,1) on $M$. If $E$ has
components $E_i^{\, h}$ and $E_i^{*h}$ relative to local coordinates $x$
and $x^*$ repectively, then

$$
p_a^*E_{i}^{*\,a} = p_aE_j^{\, b}\frac{\partial x^j}{\partial x^{*i}}.
$$
If we interpret a change of coordinates as a diffeomorphism $x^* =
x^*(x) = \varphi(x)$ we denote by $E^*$ the direct image of the tensor $E$
under the action of $\varphi$. In the case where $E$ is an almost
complex structure (that is $E^2 = -Id$), then $\varphi$ is a biholomorphism
between $(M,E)$ and $(M,E^*)$. Any (1,1) tensor field $E$ on $M$
canonically defines a contact form on $E^*M$ via

$$
\sigma = p_aE_b^{\, a}dx^b.
$$
Since 

$$
(\varphi^*)^*( p_a^*E_b^{{*}\,a}dx^{*b}) = \sigma,
$$
$\sigma$ does not depend on a choice of local coordinates (here
$\varphi^*$ is the cotangent lift of $\varphi$). Then this canonically
defines the symplectic form

$$
d\sigma = p_a\frac{\partial E_b^{\, a}}{\partial x^c}dx^c \wedge dx^b
+ E_b^{\, a}dp_a \wedge dx^b.
$$
The cotangent lift $\varphi^*$ of a diffeomorphism $\varphi$ is a
symplectomorphism for $d\sigma$. We may write

$$d\sigma = (1/2)\tau_{CB}dx^C \wedge dx^B$$
where $x^{\bar i} = p_i$; so we have

$$
\tau_{ji} = p_a \left ( \frac{\partial E_i^{\, a}}{\partial x^j} -
\frac{\partial E_j^{\, a}}{\partial x^i} \right )
$$

$$
\tau_{\bar{j} i} = E_i^{\, j}
$$

$$
\tau_{j\bar{i}} = -E_j^{\, i}
$$

$$
\tau_{\bar{j}\bar{i}} = 0
$$

We write $\widehat{E}$ for the tensor field of type (1,1) on $T^*(M)$ whose
components $\widehat{E}_B^{{}A}$  are given by 

$$\widehat{E}_B^{{}A} = \tau_{BC}\varepsilon^{CA}.$$

Thus 
 
$$
\widehat{E}_i^{\, h} = E_i^{\, h}, \ \widehat{E}_{\bar{i}}^h = 0
$$
and
$$
\widehat{E}_i^{\, \bar{h}} = 
p_a \left( \frac{\partial E_i^{\, a}}{\partial x^j} -
\frac{\partial E_j^{\, a}}{\partial x^i} \right ), 
\widehat{E}_{\bar{i}}^{{}\bar{h}} = E_h^{{}i}.
$$
In the matrix form we have

$$
\widehat{E} = \left(
\begin{matrix}
E_i^{\, h} & 0\\
 p_a \left ( \frac{\partial E_i^{\, a}}{\partial x^j} -
\frac{\partial E_j^{\, a}}{\partial x^i} \right ) & E_h^{\, i}
\end{matrix}
\right).
$$

By construction, the complete lift $\widehat{E}$ has the following
{\it invariance property}: if $\varphi$ is a local diffeomorphism of $M$
transforming $E$ to $E'$, then the direct image of $\widehat{E}$ under the
cotangent lift $\psi: = \varphi^*$ is $\widehat{E'}$.

We point out that in general, $\widehat{E}$ is not an almost complex
structure, even if $E$ is. However, the invariance condition

\begin{equation}\label{belt}
d\varphi^* \circ \widehat{E} = \widehat{E'}(\varphi^*) \circ d\varphi^*
\end{equation}
gives a Beltrami-type equation quite similarly to the almost complex case. 

\subsection{Structure properties of the Riemann map}
Assume now that $M \subset \mathbb C^n$ and let $J_\lambda$ be an almost 
complex deformation of the standard structure on $M$. 
We denote by $T_\lambda$ the elliptic prolongation $\widehat{J_\lambda}$ of
$J_{\lambda}$ ($\widehat{J_\lambda}$ is the vertical lift of $J_\lambda$ to
the cotangent bundle). So the Beltrami-type equation~(\ref{belt}) is just the
$\bar\partial_{T_\lambda}$-equation. In order to define a
biholomorphically invariant family of corresponding canonical discs we
need to consider an invariant boundary problem for the operator
$\bar\partial_{T_{\lambda}}$. 

We point out that the notion of the conormal bundle can
be easily carried to the case of an almost complex manifold. 
Let $i: T^*_{(1,0)}(M,J)\rightarrow T^*(M)$ be the canonical
identification.  In the canonical complex coordinates $(z,t)$ on
$T^*_{(1,0)}(M,J)$ an element of the fiber over the point $z$ can be
written in the form $\sum_j t_j dz^j$.  Let $D$ be a smoothly bounded
domain in $M$ with the boundary $\Gamma$. The conormal bundle
$\Sigma_J(\Gamma)$ of $\Gamma$ is a real subbundle of
$T^*_{(1,0)}(M,J)|_\Gamma$ whose fiber at $z\in \Gamma$ is defined by
$\Sigma_{z}(\Gamma) = \{ \phi \in T^*_{(1,0)}(M,J): Re \,\phi \vert
H_{(1,0)}^J(\Gamma) = 0 \}$.  Since the form $\partial \rho$ forms a
basis in $\Sigma_{z}(\Gamma)$, every $\phi \in \Sigma_J(\Gamma)$ has
the form $\phi = c \partial \rho$, $c \in \R$.

\begin{definition}
A continuous map $f:\bar\Delta \rightarrow (\bar
D,J)$, $(J_0,J)$-holomorphic on $\Delta$, is called
a {\it stationary} disc for $(D,J)$ (or for $(\Gamma,J)$) if there
exists a smooth map $\hat{f}= (f,g): \Delta
\rightarrow T^*_{(1,0)}(M,J)$, $\hat{f} \neq 0$ which is continuous on
$\bar \Delta$ and such that
\begin{itemize}
\item[(i)] $\zeta \mapsto \hat{f}(\zeta)$ satisfies the
$\bar\partial_{T_\lambda}$-equation on $\Delta$,
\item[(ii)] $(i \circ (f,\zeta^{-1}g))(\partial \Delta) \subset
\Sigma_J(\Gamma)$.
\end{itemize}
\end{definition}

We call $\hat{f}$ a {\it lift} of $f$ to the conormal bundle of
$\Gamma$. Clearly, in view of our choice of $T_{\lambda}$ the notion
of a stationary disc is {\it invariant} in the following sense: if
$\phi$ is a $\mathcal C^1$ diffeomorphism between $\bar{D}$ and 
$\bar{D}'$ and a $(J,J')$-biholomorphism from
$D$ to $D'$, then for every stationary disc $f$ in $(D,J)$ the
composition $\phi \circ f$ is a stationary discs in $(D',J')$.

Let now $D$ coincide with the unit ball $ {\mathbb B_n}$ equipped with
an almost complex deformation $J_{\lambda}$ of the standard
structure. Then it follows by definition that stationary discs in
$({\mathbb B_n},J_{\lambda})$ may be described as solutions of a
nonlinear boundary problem $({\mathcal BP}_{\lambda})$ associated with
$T_\lambda$. The above techniques give the existence and efficient
parametrization of the variety of
stationary discs in  $({\mathbb B_n},J_{\lambda})$ for $\lambda$ small
enough. This allows to apply the definition of the Riemann map and
gives its existence. We sum up our considerations in the following

\begin{theorem}\label{TH1}
Let $J_\lambda$, $J'_{\lambda}$ be almost complex perturbations 
of the standard 
structure on $\bar{\mathbb B}_n$ and $T_\lambda$, $T'_\lambda$ be their
vertical lift to the cotangent bundle. 
The Riemann map $\Psi_{T_\lambda}$ exists 
for sufficiently small $\lambda$ and
satisfies the following properties:
\begin{itemize}
\item[(a)] $\Psi_{T_\lambda} : \bar{\mathbb B}_n \backslash \{ 0 \} 
\rightarrow \bar
\Omega^\lambda \backslash \{ 0 \}$ is a smooth diffeomorphism
\item[(b)] The restriction of $\Psi_{T_\lambda}$ on every stationary disc
through the origin is $(J_0,J_0)$ holomorphic (and even linear)
\item[(c)] $\Psi_{T_\lambda}$ commutes with biholomorphisms. More
precisely for sufficiently small $\lambda'$ and for every
$\CC^1$ diffeomorphism $\varphi$ of $\bar{\mathbb B}_n$, 
$(J_\lambda,J_{\lambda'})$-holomorphic in $\mathbb B_n$
and satisfying $\varphi(0) = 0$, we have
$$
\varphi = ( \Psi_{T'_{\lambda}})^{-1} \circ d\varphi_0 \circ \Psi_{T_\lambda}.
$$
\end{itemize}
\end{theorem}

\noindent{\it Proof of Theorem~\ref{TH1}.} Conditions (a) and (b) are 
conditions (i) and (ii) of Proposition~\ref{PR1}. 

\noindent For condition~(c), let $\varphi: ({\mathbb B_n},J) \rightarrow
({\mathbb B_n},J')$ be a $(J,J')$-biholomorphism of class $\CC^1$ on
$\bar{\mathbb B}_n$ satisfying $\varphi(0) = 0$. We know that a disc
$f_v^J$ is a canonical disc for the almost complex structure $J$ if
and only if $\varphi(f_v^J)$ is a canonical disc for the almost
complex structure $J'$. Since $\varphi(f_v^J) =
f_{d\varphi_0(v)}^{J'}$ by definition, $\Psi_J(f_v^J)(\zeta) = \zeta
\,v$ and $\Psi_{J'}(f_{d\varphi_0(v)})(\zeta)=\zeta \,d\varphi_0(v)$ by
Proposition~\ref{PR1} (iii), condition (c) follows from the following
diagram:

\bigskip
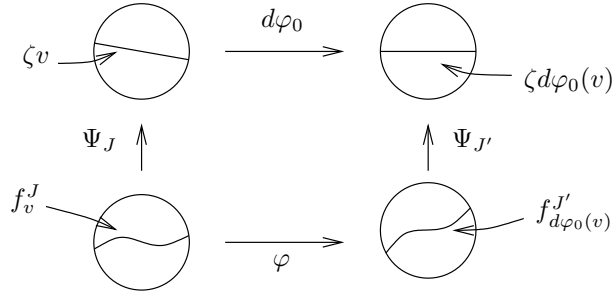
\centerline{
\begin{picture}(0,0)%
\includegraphics{diagram3.pstex}%
\end{picture}%
\setlength{\unitlength}{3947sp}%
\begingroup\makeatletter\ifx\SetFigFont\undefined%
\gdef\SetFigFont#1#2#3#4#5{%
  \reset@font\fontsize{#1}{#2pt}%
  \fontfamily{#3}\fontseries{#4}\fontshape{#5}%
  \selectfont}%
\fi\endgroup%
\begin{picture}(3300,1814)(3076,-3968)
\put(3526,-3061){\makebox(0,0)[lb]{\smash{\SetFigFont{10}{12.0}{\familydefault}{\mddefault}{\updefault}{$\Psi_J$}%
}}}
\put(3076,-3436){\makebox(0,0)[lb]{\smash{\SetFigFont{10}{12.0}{\familydefault}{\mddefault}{\updefault}{$f_v^ J$}%
}}}
\put(4651,-2311){\makebox(0,0)[lb]{\smash{\SetFigFont{10}{12.0}{\familydefault}{\mddefault}{\updefault}{$d\varphi_0$}%
}}}
\put(5851,-3061){\makebox(0,0)[lb]{\smash{\SetFigFont{10}{12.0}{\familydefault}{\mddefault}{\updefault}{$\Psi_{J'}$}%
}}}
\put(6376,-3511){\makebox(0,0)[lb]{\smash{\SetFigFont{10}{12.0}{\familydefault}{\mddefault}{\updefault}{$f_{d\varphi_0(v)}^ {J'}$}%
}}}
\put(4726,-3811){\makebox(0,0)[lb]{\smash{\SetFigFont{10}{12.0}{\familydefault}{\mddefault}{\updefault}{$\varphi$}%
}}}
\put(6301,-2686){\makebox(0,0)[lb]{\smash{\SetFigFont{10}{12.0}{\familydefault}{\mddefault}{\updefault}{$\zeta d{\varphi_0(v)}$}%
}}}
\put(3151,-2536){\makebox(0,0)[lb]{\smash{\SetFigFont{10}{12.0}{\familydefault}{\mddefault}{\updefault}{$\zeta v$}%
}}}
\end{picture}

}

Theorem~\ref{TH1} gives the main structure properties of the Riemann
map.

\subsection{Regularity of diffeomorphisms}
Riemann maps are useful for the boundary study of biholomorphisms in
almost complex manifolds. We have
\begin{corollary}
If $\lambda <<1$ and $\varphi$ is a $\CC^1$ diffeomorphism of 
$\bar{\mathbb B}_n$, 
$(J_\lambda,J_{\lambda'})$-holomorphic in $\mathbb B_n$
satisfying $\varphi(0) = 0$, then $\varphi$ is of class $\CC^{\infty}$
on $\bar{\mathbb B}_n$.
\end{corollary}
\proof This follows immediately by Theorem~\ref{TH1} condition~(c) 
since the Riemann map is smooth up to the boundary. \qed

\vskip 0,1cm
Every almost complex structure is locally a small deformation 
of the standard structure. So we have the following partial generalization of
Fefferman's theorem:

\begin{theorem}\label{CO1}
Let $(M,J)$, $(M,J')$ be almost complex manifolds and let $D$ and
$D'$ be domains in $M$ and $M'$ respectively. Assume that $\partial D$ and
$\partial D'$ contain smooth strictly pseudoconvex open pieces $\Gamma$
and $\Gamma'$. If $\varphi$ is a $\mathcal C^1$ diffeomorphism between
$D \cup \Gamma$ and $D' \cup \Gamma'$, 
$(J,J')$-holomorphic on $D$ such that
$\varphi(\Gamma) \subset \Gamma'$, 
then $\varphi$ is of class $\CC^{\infty}$ on $D \cup \Gamma$.
\end{theorem}

\proof Let $q$ be a point of $\Gamma$ and $q' = \phi(q)$. We may
assume that $J(q) = J_0$ and $J'(q') = J_0$. Since $\Gamma$ and
$\Gamma'$ are strictly pseudoconvex, in local coordinates near $q$
(resp. $q'$) satisfying $z(q) = 0$ (resp. $z(q') = 0$) they are small
deformations of the Siegel sphere $\mathbb H = \{ Re z^n + \parallel z
\parallel^2 = 0 \}$. After the Caley map, $\mathbb H$ is transformed
to the unit sphere $\mathbb S$ and $\Gamma$ (resp. $\Gamma'$) may be
viewed as a small deformation of the sphere along the boundary of some
stationary disc $f^0$, $f^0 (0) = 0$ such that the corresponding disc
on the Siegel domain $\{ Re z^n + \parallel z \parallel^2 < 0 \}$ has
direction at the center parallel to the holomorphic tangent space of
$\mathbb H$ at the origin.  Denoting again by $J$ resp. $J'$ the
images of our almost complex structures under the Caley map, we see
that in a neighborhood of the disc $f^0$ they may be viewed as an
arbitrarily small deformations of the standard structure. So we may
consider as in Subsection~4.4 the Riemann map $\Psi_{T_\lambda,v^0}$ 
where $v^0$ is the tangent vector at the origin of the disc $f^0$.
Since the property (c) of Theorem~\ref{TH1} obviously holds
for local Riemann maps, we conclude. \qed

\vskip 0,1cm
Theorem~\ref{CO1} admits the following global version.  

\begin{theorem}\label{CO2}
Let $(M,J)$ and $(M,J')$ be almost complex manifolds. Let $D$ and
$D'$ be smoothly bounded strictly pseudoconvex domains in $M$ and $M'$
respectively. Every $\CC^1$ diffeomorphism between $\bar{D}$ and $\bar{D}'$,
$(J,J')$-holomorphic on $D$, is of class $\CC^{\infty}$ on $\bar D$.
\end{theorem}

The following reformulation of Theorem~\ref{CO1} may be considered as a
geometrical version of the elliptic regularity for manifolds with boundary.

\begin{theorem}
Let $M$ and $M'$ be two $\mathcal C^{\infty}$ smooth real
$2n$-dimensional manifolds, $D \subset M$ and $D' \subset M'$ be
relatively compact domains. Suppose that there exists an almost
complex structure $J$ of class $\mathcal C^{\infty}$ on $\bar D$ such
that $(D,J)$ is strictly pseudoconvex. Then a $\mathcal C^{1}$
diffeomorphism $\phi$ between $\bar{D}$ and $\bar{D}'$ is of class 
$\mathcal C^{\infty}(\bar D)$ if and only if the
direct image $J':= \phi^*(J)$ is of class $\mathcal
C^{\infty}(\bar{D'})$ and $(D',J')$ is strictly pseudoconvex.
\end{theorem}

\subsection{Rigidity and local equivalence problem}

Condition $(c)$ of Theorem~\ref{TH1} implies the following partial 
generalization of Cartan's theorem for almost complex manifolds:
\begin{corollary}
If $\lambda <<1$ and if 
$\varphi$ is a $\CC^1$ diffeomorphism of $\bar{\mathbb B}_n$, 
$(J_\lambda,J_{\lambda'})$-holomorphic in $\mathbb B_n$,
satisfying $\varphi(0) = 0$ and $d\varphi(0) = I$ then $\varphi$ is the 
identity.
\end{corollary}
This provides an efficient parametrization of the
isotropy group of the group of biholomorphisms of $({\mathbb B_n},J_\lambda)$.

\vskip 0,1cm
We can solve the local biholomorphic 
equivalence problem between almost complex manifolds in terms of the 
Riemann map similarly to~\cite{bl-du-ka87,le88} (see the paper
\cite{li50} by P. Libermann for a traditional approach to this problem 
based on Cartan's equivalence method for $G$-structures).
Let $I^\lambda$ (resp. $(I')^\lambda$) be the indicatrix of 
($\mathbb B_n,J_\lambda$) (resp. ($\mathbb B_n,J'_\lambda$)) bounding
the domain $\Omega^\lambda$ (resp. $(\Omega')^\lambda$) and let 
$\Psi_{T_\lambda}$ (resp. $\Psi_{T'_\lambda}$) be the associated Riemann map.
This induces the almost complex structure 
$J_\lambda^*:=d\Psi_{T_\lambda} \circ J_\lambda \circ d(\Psi_{T_\lambda})^{-1}$
 (resp.
$(J'_\lambda)^*:= 
d\Psi_{T'_\lambda} \circ J_\lambda \circ d(\Psi_{T'_\lambda})^{-1}$) 
on $\Omega^{\lambda}$ (resp. $(\Omega')^{\lambda}$). Then we have:

\begin{theorem}\label{TH3}
The following conditions are equivalent:

$(i)$ There exists a $\CC^\infty$ diffeomorphism $\varphi$ of 
$\bar{\mathbb B}_n$, $(J_\lambda,J_{\lambda}')$-holomorphic on $\mathbb B_n$
and satisfying $\varphi(0)=0$, 

$(ii)$ There exists a $J_0$-linear isomorphism
$L$ of $\C^n$, $(J_\lambda^*,(J'_\lambda)^*)$-holomorphic 
on $\Omega^\lambda$ and such that $L(\Omega^\lambda)=(\Omega')^\lambda$.
\end{theorem}

\noindent{\it Proof}. If $\varphi$ satisfies condition $(i)$, 
the commutativity of the following diagram (in view of
Theorem~\ref{TH1})

\bigskip
\centerline{
\begin{picture}(0,0)%
\includegraphics{commut3.pstex}%
\end{picture}%
\setlength{\unitlength}{3947sp}%
\begingroup\makeatletter\ifx\SetFigFont\undefined%
\gdef\SetFigFont#1#2#3#4#5{%
  \reset@font\fontsize{#1}{#2pt}%
  \fontfamily{#3}\fontseries{#4}\fontshape{#5}%
  \selectfont}%
\fi\endgroup%
\begin{picture}(2898,1857)(4340,-3863)
\put(7238,-3034){\makebox(0,0)[lb]{\smash{\SetFigFont{9}{10.8}{\familydefault}{\mddefault}{\updefault}{$\Psi'_{T_\lambda}$}%
}}}
\put(6840,-3829){\makebox(0,0)[lb]{\smash{\SetFigFont{9}{10.8}{\familydefault}{\mddefault}{\updefault}{$(\mathbb B_n,J'_\lambda)$}%
}}}
\put(5704,-3545){\makebox(0,0)[lb]{\smash{\SetFigFont{9}{10.8}{\familydefault}{\mddefault}{\updefault}{$\varphi$}%
}}}
\put(4340,-3829){\makebox(0,0)[lb]{\smash{\SetFigFont{9}{10.8}{\familydefault}{\mddefault}{\updefault}{$(\mathbb B_n,J_\lambda)$}%
}}}
\put(4795,-3034){\makebox(0,0)[lb]{\smash{\SetFigFont{9}{10.8}{\familydefault}{\mddefault}{\updefault}{$\Psi_{T_\lambda}$}%
}}}
\put(4340,-2295){\makebox(0,0)[lb]{\smash{\SetFigFont{9}{10.8}{\familydefault}{\mddefault}{\updefault}{$(\Omega^ \lambda,J_\lambda^ *)$}%
}}}
\put(6897,-2295){\makebox(0,0)[lb]{\smash{\SetFigFont{9}{10.8}{\familydefault}{\mddefault}{\updefault}{$((\Omega')^ \lambda,(J')_\lambda^ *)$}%
}}}
\put(5533,-2125){\makebox(0,0)[lb]{\smash{\SetFigFont{9}{10.8}{\familydefault}{\mddefault}{\updefault}{$L=d\varphi_0$}%
}}}
\end{picture}

}
\bigskip
\noindent shows that $L:=d\varphi_0$ satisfies condition $(ii)$. 
Conversely if $L$ satisfies condition $(ii)$ 
then the map $\varphi:=(\Psi'_{T_\lambda})^{-1} \circ
L \circ \Psi_{T_\lambda}$ satisfies condition $(i)$. \qed

\vskip 0,3cm
In conclusion we point out that there are many open questions concerning
the Riemann map on almost complex manifolds (contact properties, relation
with the Monge-Amp\`ere equation, \dots). 
They will be studied in a forthcoming paper.

\end{document}